\documentclass[11pt]{amsart}

\usepackage[centering]{geometry}            

\usepackage{graphicx}
\usepackage{listings}
\usepackage{amssymb}
\usepackage{dsfont}
\usepackage{url}
\usepackage{verbatim}

\usepackage{amssymb,amsmath,setspace,geometry,indentfirst,changepage,inputenc}
\usepackage{mathrsfs,amsfonts,savesym,graphicx,bm,amsthm}

\usepackage[all]{xy}
\usepackage{mathrsfs}
\usepackage{color}

\usepackage{dcolumn}
\usepackage{cite}

\usepackage{booktabs}
\usepackage{diagbox} 
\usepackage{multirow} 
\usepackage{makecell} 
\usepackage{longtable} 
\usepackage{array}
\usepackage{tabularx}
\usepackage{color}

\usepackage{amsmath,amsthm}

\usepackage[mathscr]{eucal}
\usepackage[all]{xy}
\usepackage{mathrsfs}
\usepackage{color}

\newtheorem{theorem}{Theorem}[section]
\newtheorem{lemma}[theorem]{Lemma}

\newtheorem{conjecture}[theorem]{Conjecture}
\newtheorem{corollary}[theorem]{Corollary}
\newtheorem{proposition}[theorem]{Proposition}

\newtheorem{question}[theorem]{Question}
\newtheorem{problem}[theorem]{Problem}
\newtheorem{definition-lemma}[theorem]{Definition-Lemma}
\newtheorem{definition-theorem}[theorem]{Definition-Theorem}

\theoremstyle{definition}
\newtheorem{example}[theorem]{Example}
\newtheorem{definition}[theorem]{Definition}

\newtheorem{remark}[theorem]{Remark}

\newtheorem*{ack}{Acknowledgements}

\numberwithin{equation}{section}

\allowdisplaybreaks

\title{On the Tensor Property of Bernstein-Sato Polynomials}

\author{Quan Shi}
\address{
		Department of Mathematical Sciences,
		Tsinghua University,
		Beijing, 100084, P. R. China.}
\email{shiq24@mails.tsinghua.edu.cn}

\author{Huaiqing Zuo}
\address{Department of Mathematical Sciences,
	Tsinghua University,
	Beijing, 100084, P. R. China.}
\email{hqzuo@mail.tsinghua.edu.cn}

\begin{document}

	\begin{abstract}
		We prove the multiplicative Thom–Sebastiani rule for Bernstein–Sato polynomials, \textcolor{black}{thereby} answering a long-standing question of Budur and Popa. 
		\textcolor{black}{This result is extended to the tensor product} of two effective divisors on the product of two arbitrary nonsingular complex varieties. \textcolor{black}{As a consequence, we also derive a multiplicative Thom-Sebastiani type result \textcolor{black}{for} the strong monodromy conjecture. Furthermore, we propose analogous conjectures for Bernstein–Sato polynomials for ideals and verify them in the case of monomial ideals.}
		
		
	\end{abstract}
	\maketitle
	
	\tableofcontents
	
	\subsection*{Convention} In this paper, we adopt \textcolor{black}{the} multi-index \textcolor{black}{notation}. We use $\bm x$ to denote an array of variables or numbers depending on the context. \textcolor{black}{For a tuple $\bm x = (x_1,...,x_n)$ and $\bm \alpha = (\alpha_1,...,\alpha_n) \in \mathbb Z^n$, we write $\bm x^{\bm \alpha}:=\prod_i x_i^{\alpha_i}$ and $\vert \bm x\vert := x_1+...+x_n$.} 
	
	\section{Introduction}\label{introduction}
	
	Suppose $X$ is a nonsingular \textcolor{black}{integral} complex variety and $f$ is a non-zero regular function. \textcolor{black}{Let $\mathscr D_X$ be the sheaf of differential operators on $X$. This is the sheaf of $\mathbb C$-subalgebras of $\mathcal End_{\mathbb C}(\mathcal O_X)$ generated by $\mathcal O_X$ and $\mathcal T_X$, the tangent sheaf of $X$.} The Bernstein-Sato polynomial $b_f(s)$ is the non-zero \textcolor{black}{monic} polynomial in $s$ of the minimal degree such that
	\begin{equation}
		b_f(s) \cdot f^s \in \mathscr D_X[s] \cdot f^{s+1}, \label{eq1}
	\end{equation}
	where $s$ is a \textcolor{black}{variable, $f^s$ is a formal symbol,} and $\mathscr D_X[s]$ \textcolor{black}{acts} on the right-hand side in the expected way. For $X = \mathbb A_{\mathbb C}^{n}$, the condition \eqref{eq1} can be rewritten as 
	\begin{displaymath}
		b_f(s) \cdot f^s =\textcolor{black}{ P(s,x_1,...,x_n,\partial_{x_1},...,\partial_{x_n})\cdot f^{s+1},}
	\end{displaymath}
	where $P$ is an element in \textcolor{black}{$\Gamma(\mathbb A_{\mathbb C}^n,\mathscr D_{\mathbb A_{\mathbb C}^n})[s] = A_n(\mathbb C)[s]$, the Weyl algebra tensored with $\mathbb C[s]$.} The existence of $b_f(s)$ is a non-trivial result, which follows from the holonomicity of \textcolor{black}{the} $\mathscr D_{X_s}$-module $\mathcal O_{X_s}[\frac{1}{f}] \cdot f^s$, where $X_s = X\times_{\mathrm{Spec}\, \mathbb C} \mathrm{Spec}\, \mathbb C(s)$, \textcolor{black}{see subsection \ref{subSec Bernstein-Sato polynomial, Existence and Some Results}}.

	
	\textcolor{black}{The Bernstein-Sato polynomial captures profound geometric data \textcolor{black}{regarding} hypersurface singularities. For example, it is the algebraic avatar of the Milnor fiber. This is established by the well-known Kashiwara-Malgrange theorem, saying that, after applying the transformation $\lambda \mapsto e^{2\pi \sqrt{-1}\lambda}$, the roots of $b_f(s)$ correspond precisely to monodromy eigenvalues of cohomologies of Milnor fibers, see \cite{Kashiwara_Malgrange_Theorem_Kashiwara,Kashiwara_Malgrange_Theorem_Malgrange}. Additionally, the Bernstein-Sato polynomial serves as a key invariant in the classification of singularities, see \cite{KS11,Sai93,Sai09}.}

	The definition \textcolor{black}{of the Bernstein-Sato polynomial} can be extended \textcolor{black}{globally to a non-zero effective Cartier divisor $D$ on $X$.} \textcolor{black}{Choose} an arbitrary affine covering $\{U_i\}_i$ and sections $f_i \in \Gamma(U_i,\mathcal O_X)$ such that $D|_{U_i} = \mathrm{div}(f_i)$, then one defines $b_D(s) := \mathrm{lcm}_i\, b_{f_i}(s)$. \textcolor{black}{The global Bernstein-Sato polynomial $b_D(s)$ is deeply connected to other fundamental objects in singularity theory and birational geometry, e.g. the multiplier ideal, see \cite{ELSV04}. }
	
	\vspace{0.5em}

	For \textcolor{black}{non-zero polynomials} $f\in \textcolor{black}{\mathbb C[x_1,...,x_n]}$ and $g \in \textcolor{black}{\mathbb C[y_1,...,y_m]}$, Nero Budur proposed the following problem.
	\begin{problem}[{\cite[2.11]{Bud12}}]\label{Budur's conjecture}
		Is there a multiplicative Thom-Sebastiani rule for Bernstein-Sato polynomials? That is, is there any relation between $b_{f}(s)$, $b_g(s)$, and $b_{f\cdot g}(s)$?
	\end{problem}
	Years later, Mihnea Popa further asked the following. 
	\begin{problem}[{\cite[page 33]{DMBG}}]\label{Popa' conjecture}
		Do the Bernstein-Sato polynomials of the above $f$, $g$, and $f\cdot g$ satisfy
		\begin{equation}\label{eq2}
			b_{f}(s)\cdot b_g(s) = b_{f\cdot g}(s) ?
		\end{equation}
	\end{problem} 
	\textcolor{black}{Although one divisibility, namely $b_{f\cdot g}(s) \mid b_{f}(s)\cdot b_g(s)$, is straightforward, the converse direction remains highly non-trivial.} In literature, only few clues related to \textcolor{black}{Problem \ref{Budur's conjecture} and Problem \ref{Popa' conjecture}} can be found. 
	\textcolor{black}{\cite[Proposition 2.6.1]{Walter_PhD} implies} that $b_f(s)\cdot b_g(s) = b_{f\cdot g}(s)$ holds for $f$ and $g$ whose Bernstein-Sato polynomials can be achieved by differential operators involving only $s$ and \textcolor{black}{ ${\partial_{x_1}},...,\partial_{x_n},{\partial_{y_1}},...,\partial_{y_m}$}. In \cite{DMBG}, \textcolor{black}{the formula for the Bernstein-Sato polynomial of a monomial helps prove the identity when either $f$ or $g$ is a monomial. }
	
	The problems are related to the tensor property of Bernstein-Sato polynomials since $f\cdot g \in \textcolor{black}{\mathbb C[x_1,...,x_n,y_1,...,y_m]}$ can be identified with $f\otimes g \in \textcolor{black}{\mathbb C[x_1,...,x_n] \otimes_{\mathbb C} \mathbb C[y_1,...,y_m]}$. \textcolor{black}{In this paper, we prove the following more general result.}
	\textcolor{black}{\begin{theorem}\label{tensor theorem for bernstein-sato polynomial}
			Suppose $(X_1,D_1)$ and $(X_2,D_2)$ are pairs consisting of a nonsingular integral complex variety and a non-zero effective Cartier divisor. Let $D_1 \otimes D_2$ be the Cartier divisor on $X_1 \times X_2$ locally defined by the tensor product of local equations of $D_1$ and $D_2$. Then we have
			\begin{displaymath}
				b_{D_1}(s) \cdot b_{D_2}(s) = b_{D_1\otimes D_2} (s).
			\end{displaymath}
	\end{theorem}}
	\textcolor{black}{Setting $(X_1,D_1) = (\mathbb A_{\mathbb C}^n,\mathrm{div}(f))$ and $(X_2,D_2) = (\mathbb A_{\mathbb C}^m,\mathrm{div}(g))$, we address Problem \ref{Budur's conjecture} and Problem \ref{Popa' conjecture}.
		\begin{corollary}
			For non-zero polynomials $f\in \mathbb C[x_1,...,x_n]$ and $g \in \mathbb C[y_1,...,y_m]$, we have $b_f(s) \cdot b_g(s) = b_{f\cdot g}(s)$.
	\end{corollary}}

	%
	%
	%
	%

	%

	\vspace{0.5em}
	
	
	\textcolor{black}{Next, we present an application of Theorem \ref{tensor theorem for bernstein-sato polynomial} to the monodromy conjecture, one of the main open problems in singularity theory. The following are its $p$-adic and motivic version by Igusa and Denef-Loeser respectively, see \cite{Nic10,Vey24} for more details.}

	Let $p$ be a prime number. The $p$-adic number field $\mathbb Q_{p}$ is a locally compact group and admits a Haar measure $\mu_p$, \textcolor{black}{unique up to scalar, see \cite{MR0234930}}. After normalizing, one has $\mu_p(p^r\mathbb Z_p) = p^{-r}$ for all $r\in \mathbb Z$. \textcolor{black}{The norm associated to $\mu_p$ satisfies $\vert p^{r}\vert_p = p^{-r}$. }
	
	\textcolor{black}{Suppose $f\in \mathbb Z[x_1,...,x_n]$ is a non-zero polynomial. One can define a function of $s$, holomorphic on $\{s \in \mathbb C \mid \mathrm{Re}\, s > 0\}$, as follows.}
	\begin{displaymath}
		Z_p(f;s) := \int_{\mathbb Z_p^n} \vert f\vert_p^s \, d\mu_p.
	\end{displaymath}
	Igusa showed that $Z_p(f;s)$ \textcolor{black}{admits a meromorphic extension to the whole complex plane and} is a rational function of $p^{-s}$\textcolor{black}{, cf. \cite{Igu75}}. \textcolor{black}{The extended meromorphic function, also denoted by $Z_p(f;s)$, is called the $p$-adic zeta function of $f$. }
	\begin{conjecture}[$p$-adic strong monodromy conjecture, \cite{Igu00}]\label{StrMonoConj}
		Given a non-zero polynomial $f\in \mathbb Z[x_1,...,x_n]$, for every prime $p$ large enough, if $s_0$ is a pole of $Z_p(f;s)$, then \textcolor{black}{$\mathrm{Re}\,{s_0}$} is a root of $b_f(s)$. Moreover, if the order of $s_0$ as a pole is $m$, then $\textcolor{black}{\mathrm{Re}\,{s_0}}$ is a root of $b_f(s)$ of multiplicity $\geq m$.
	\end{conjecture} 
	
	\textcolor{black}{Let $h\in \mathbb C[x_1,...,x_n]$ be a non-zero polynomial. The motivic zeta function of $h$, denoted by $Z^{\mathrm{mot}}(h;s)$, is an analogue of $Z_p(f;s)$. Intrinsically, $Z^{\mathrm{mot}}(h;s)$ is defined as the normalized generating series of classes of contact loci in the Grothendieck ring, but here for simplicity, we only provide the following extrinsic formula, cf. \cite{DL98}.
		\begin{theorem}[Denef-Loeser] 
			Notations are as above. Suppose $\mu : X\to \mathbb C^n$ is a log resolution of $(\mathbb C^n,\mathrm{div}(h))$ and $\mu^*\mathrm{div}(h)  = \sum_{i\in S} N_i E_i,
			K_{\mu}  = \sum_{i\in S} (b_i-1) E_i,$
			where $\{E_i\}_{i\in S}$ have simple normal crossings and $K_{\mu}$ is the relative canonical divisor. Then
			\begin{displaymath}
				Z^{\mathrm{mot}}(h;s) = \mathbb L^{-n} \sum_{I\subseteq S} [E_I^\circ] \cdot \prod_{i\in I}\frac{(\mathbb L-1)\mathbb L^{-(N_is+b_i)}}{1-\mathbb L^{-(N_is+b_i)}}.
			\end{displaymath}
			Here, $E_I^\circ = (\bigcap_{i\in I} E_i) \setminus (\bigcup_{j\in S\setminus I} E_j)$, $[-]$ means the class in the Grothendieck ring of complex varieties, and $\mathbb L := [\mathbb A^1_{\mathbb C}]$.
		\end{theorem}
		\begin{conjecture}[motivic strong monodromy conjecture, \cite{DL98}]\label{StrMonoConj_Denef_Loeser}
			Given a non-zero polynomial $h\in \mathbb C[x_1,...,x_n]$, if $s_0$ is a pole of $Z^{\mathrm{mot}}(h;s)$, then ${s_0}$ is a root of $b_h(s)$. Moreover, if the order of $s_0$ as a pole is $m$, then ${s_0}$ is a root of $b_h(s)$ of multiplicity $\geq m$.
		\end{conjecture} 
		Both Conjecture \ref{StrMonoConj} and Conjecture \ref{StrMonoConj_Denef_Loeser} remain widely open. One major obstacle lies in the insufficient understanding of the cancellation of poles of these zeta functions.}
	

	\textcolor{black}{Theorem \ref{tensor theorem for bernstein-sato polynomial} immediately yields the following multiplicative Thom-Sebastiani type result for these strong monodromy conjectures.
		\begin{corollary}\label{product monodromy conjecture}
			Let $\mathbb K = \mathbb Z$ (resp. $\mathbb C$), and let $f \in \mathbb K[x_1, ..., x_n]$ and $g \in \mathbb K[y_1, ..., y_m]$ be non-zero polynomials. If Conjecture \ref{StrMonoConj} (resp. Conjecture \ref{StrMonoConj_Denef_Loeser}) holds for $f$ and $g$, then the corresponding conjecture also holds for $f\cdot g \in \mathbb K[x_1,....x_n,y_1,...,y_m]$.
	\end{corollary}}

	
	\vspace{0.5em}

	The Bernstein-Sato polynomial can be generalized to \textcolor{black}{the case of ideals, see \cite{BMS06b,Mus22}}. Let $\mathfrak a \subseteq \mathbb C[x_1,...,x_n]$ \textcolor{black}{be a non-zero ideal with non-zero generators $f_1,...,f_r$}, then \textcolor{black}{$b_{\mathfrak a}$ is the monic polynomial in one variable of the minimal degree} such that
	\begin{displaymath}
		\textcolor{black}{b_{\mathfrak a}(s_1+...+s_r) f_1^{s_1}f_2^{s_2} ... f_r^{s_r} \in \sum_{\bm u\in \mathbb Z^r, \vert \bm u\vert = 1} A_n(\mathbb C)[s_1,...,s_r] \cdot (\prod_{u_i < 0}\binom{s_i}{-u_i}) f_1^{s_1+u_1}f_2^{s_2+u_2}...f_r^{s_r+u_r}.}
	\end{displaymath}
	The general computation of $b_{\mathfrak a}$ is very difficult. \textcolor{black}{However, for monomial ideals,} we have a \textcolor{black}{combinatorial} description, \textcolor{black}{see \cite{BMS06c}. More details can be found in {subsection \ref{subSec Bernstein-Sato ideal case}}.}
	
	One may be interested in analogues of {Problems \ref{Budur's conjecture}} and {Problems \ref{Popa' conjecture}} \textcolor{black}{for ideals}. \textcolor{black}{That is, for non-zero ideals $\mathfrak a\subseteq \mathbb C[x_1,...,x_n]$ and $\mathfrak b\subseteq \mathbb C[y_1,...,y_m]$, one may ask whether there is} a relation between $b_{\mathfrak a},b_{\mathfrak b}$, and $b_{\mathfrak a\cdot \mathfrak b}$, or \textcolor{black}{even if one has $b_{\mathfrak a} \cdot b_{\mathfrak b} = b_{\mathfrak a\cdot \mathfrak b}$.} \textcolor{black}{We show that the equality holds when one of the ideals is principal.}
	\begin{proposition}\label{generalized Popa's conjecture in a case}
		Let $\mathfrak a \subseteq \textcolor{black}{\mathbb C[x_1,...,x_n]}$ be a non-zero ideal and $g \in \textcolor{black}{\mathbb C[y_1,....,y_m]}$ be a non-zero polynomial. Then $b_\mathfrak a \cdot b_g = b_{\mathfrak a \cdot (g)}$.
	\end{proposition}
	

	\textcolor{black}{The problem remains interesting even if one omits multiplicities. Let $W_{\mathfrak a},W_{\mathfrak b}$, and $W_{\mathfrak a\cdot \mathfrak b}$ denote the set of roots of $b_\mathfrak a,b_{\mathfrak b}$, and $b_{\mathfrak a\cdot \mathfrak b}$ respectively. For monomial ideals, the following Thom-Sebastiani rule holds.
		\begin{theorem}\label{bab for monomial ideals1}
			Let $\mathfrak a\subseteq \mathbb C[x_1,...,x_n]$ and $\mathfrak b\subseteq \mathbb C[y_1,....,y_m]$ be non-zero monomial ideals, then
			\begin{displaymath}
				W_{\mathfrak a}\cup W_\mathfrak b = W_{\mathfrak a \cdot \mathfrak b}.
			\end{displaymath}
	\end{theorem}}
	
	%
	%
	%
	
	\textcolor{black}{Based on Proposition \ref{generalized Popa's conjecture in a case} and Theorem \ref{bab for monomial ideals1}, we propose the following conjecture concerning Bernstein–Sato polynomials for general ideals.}
	\begin{conjecture}\label{WabWaWb conjecture}
		\textcolor{black}{Let $\mathfrak a\subseteq \mathbb C[x_1,...,x_n]$ and $\mathfrak b\subseteq \mathbb C[y_1,....,y_m]$ be non-zero ideals. Then}
		
		\noindent\textcolor{black}{(1)(weak) $W_{\mathfrak a}\cup W_\mathfrak b =  W_{\mathfrak a \cdot \mathfrak b} \mod \mathbb Z$;}
		
		\noindent\textcolor{black}{(2)(strong) $W_{\mathfrak a}\cup W_\mathfrak b = W_{\mathfrak a \cdot \mathfrak b}$.}
	\end{conjecture}

%

\textcolor{black}{The organization of this paper is as follows. {Section \ref{Preliminary}} is for preliminary material on Bernstein-Sato polynomials. In {Section \ref{On Popa's Conjecture}}, we present the motivations, proofs, and generalization of \textcolor{black}{both} {Problem \ref{Budur's conjecture}} and {Problem \ref{Popa' conjecture}}.} 

\begin{remark}
	\textcolor{black}{After this paper was finished and submitted, \cite{lee2024multiplicative} appeared. It also solves Problem \ref{Budur's conjecture} and Problem \ref{Popa' conjecture}, and has some overlaps with our paper. However, two solutions are completed independently using different methods. Our solution arises from a characterization of the Bernstein-Sato polynomial via commutative algebra.}
	%
	%
\end{remark}
\vspace{-1.8em}

\textcolor{black}{\begin{ack}
		We thank the anonymous referees for their valuable comments, which have helped us improve the quality of the paper, especially regarding Theorem \ref{bab for monomial ideals1} and Conjecture \ref{WabWaWb conjecture}.
		Q. Shi and H. Zuo were supported by BJNSF Grant 1252009.  H. Zuo was supported by NSFC Grant 12271280. We thank Siyong Tao for discussion.
\end{ack}}


\section{\textcolor{black}{Preliminaries for Bernstein-Sato polynomials}}\label{Preliminary}

\subsection{\textcolor{black}{Basics}}\label{subSec Bernstein-Sato polynomial, Existence and Some Results}

The Bernstein-Sato polynomial $b_f(s)$ is \textcolor{black}{a fundamental object} in $\mathscr D$-module theory. Historically, Joseph Bernstein and Mikio Sato discovered $b_f(s)$ independently in the 1970s. Bernstein found it when studying the distribution $f^s$, and Sato arrived at it when studying prehomogeneous vector spaces\textcolor{black}{, see \cite{doi:10.1142/S0219199721500838}. }

\textcolor{black}{The existence of $b_f(s)$ is itself non-trivial. For $X = \mathbb A_{\mathbb C}^n$, the existence is proved by using the Bernstein filtration on the Weyl algebra $A_n(\mathbb C(s))$, see \cite{MR0290097}. For a general nonsingular variety $X$, there is no such filtration. Alternatively, one should carefully establish the theory of holonomic $\mathscr D_{X_s}$-modules, where $p : X_s := X \times_{\mathrm{Spec}\, \mathbb C} \mathrm{Spec}\, \mathbb C(s) \to X$ is the base change and $\mathscr D_{X_s} := p^*\mathscr D_X$, see \cite{DMBG}.}

%
%

\textcolor{black}{Apart from monodromy eigenvalues, the Bernstein-Sato polynomial is also relevant to other important invariants of singularities, e.g. archimedean zeta functions, period integrals, higher multiplier ideals, see \cite{DLY24,Bla24,SY23} for more details.}

\textcolor{black}{The Bernstein–Sato polynomial admits several natural generalizations. First, one can replace $f$ by a tuple of polynomials $F = (f_1,...,f_r)$ and define the Bernstein-Sato ideal $B_F \subseteq \mathbb C[s_1,...,s_r]$, cf. \cite{Zero_Loci_of_Bernstein-Sato_Ideals,BSZ25}. Second, one can change $f$ to an ideal as in subsection \ref{subSec Bernstein-Sato ideal case}. Third, one can define the Bernstein-Sato polynomial for a section of $\iota_+\mathcal O_X$, where $\iota$ is the graph embedding of $f$ and $\iota_+$ is the $\mathscr D$-module push-forward functor, see \cite{DMBG}. }

\vspace{0.5em}

The computation of $b_f(s)$ \textcolor{black}{in general cases} is a rather difficult problem. \textcolor{black}{We list below some known results for special polynomials.}

\begin{theorem}[see \cite{DMBG}]
	Let $f = \det(x_{ij}) \in \mathbb C[x_{ij}]_{1\leq i,j\leq n}$ be the determinant of  \textcolor{black}{a generic matrix in $n^2$ variables}, then
	\begin{displaymath}
		b_f(s) = (s+1)(s+2)...(s+n).
	\end{displaymath}
\end{theorem}
	\begin{theorem}[\cite{MR3573952}]
		The Bernstein-Sato polynomial of a \textcolor{black}{reduced} generic central hyperplane arrangement $f \in \mathbb C[\bm x]$ of degree $l \geq n$ is
		\begin{displaymath}
			b_f(s) = (s+1)^{n-1} \prod_{j = 0}^{2l-n-2} (s+\frac{j+n}{l}).
		\end{displaymath}
	\end{theorem}
	
	\begin{theorem}[see \cite{MR2766095} or \cite{DMBG}]\label{b-function for quasi-homogeneous}
		Suppose $f\in \mathbb C[\bm x]$ is a weighted homogeneous polynomial defining an isolated singularity at the origin, with weight type $\bm w = (w_1,...,w_n)$ \textcolor{black}{such that $\bm w\cdot \bm \alpha = \sum_{i=1}^n \alpha_iw_i = 1$ for all monomials $\bm x^{\bm \alpha}$ appearing in $f$}. Let $M_f := \mathbb C\{\bm x\}/(\frac{\partial f}{\partial x_1},...,\frac{\partial f}{\partial x_n})$ be the Milnor algebra of $f$ and $\{\bm x^{\bm \alpha_i}\}_{1\leq i\leq \mu}$ be a monomial basis of $M_f$\textcolor{black}{, where $\mu = \dim_{\mathbb C} M_f$ is the Milnor number of $f$. Let $\Delta$ denote the set (without multiplicity) of all rational numbers of the form $\bm w\cdot \bm \alpha_i$.} The Bernstein-Sato polynomial of $f$ is
		\begin{displaymath}
			b_f(s) = (s+1) \prod_{\rho \in \Delta} (s+\vert \bm w\vert+\rho).
		\end{displaymath}
	\end{theorem}

	\subsection{\textcolor{black}{Characterization via commutative algebra}} \label{ch_CA} 
	
	\textcolor{black}{Let $X$ be a nonsingular integral complex variety
		and $\mathscr D_X $ be the sheaf of differential operators on $X$. Then there is a natural filtration $F_\bullet\mathscr D_X$ on $\mathscr D_X$. That is, 
		for any small open subset $U \subseteq X$ and local coordinates $x_1,...,x_n$ on it, we have $F_j\mathscr D_X|_U = \bigoplus_{\bm a \in \mathbb N^n, \vert \bm a\vert \leq j} \mathcal O_U \cdot \bm{\partial_x^{a}}$.}
	%

\textcolor{black}{Further assuming $X$ is affine, and write $A := \Gamma(X,\mathcal O_X)$ for its coordinate ring. Define $D(A) := \Gamma(X,\mathscr D_X)$ and $F_jD(A) := \Gamma(X,F_j\mathscr D_X)$. Let $f \in A$ be a non-zero regular function. In this case, the Bernstein-Sato polynomial of $f$ is the minimal monic polynomial $b_f(s)\in \mathbb C[s]$ such that there exists $P\in D(A)[s]$ satisfying}
\begin{equation}\label{eq21}
	\textcolor{black}{b_f(s)\cdot f^s = P\cdot f^{s+1}}. 
\end{equation}

\begin{remark}
	\textcolor{black}{One can define the Bernstein-Sato polynomial for a non-zero holomorphic function germ $f \in \mathbb C\{\bm x\}$ by replacing $A$ with the ring of holomorphic function germs $\mathbb C\{\bm x\}$, and $D(A)$ with $\mathbb C\{\bm x\}[\bm \partial_{\bm x}]$.} \textcolor{black}{Roots of $b_f(s)$ are also important. For example, if $f$ has an isolated singularity at the origin, then the negative of the maximal root of $b_f(s)/(s+1)$ equals the minimal spectrum number, see \cite{saito2006introduction}.}
\end{remark}

\textcolor{black}{Now, consider the $D(A)[s]$-module $D(A)[s] \cdot f^{s+1} \subseteq A_f[s]\cdot f^s$. For all $j\geq 0$, we have an ideal $I_j$ of $A[s]$ satisfying}
\begin{displaymath}
	F_jD(A)[s]  \cdot f^{s+1} = I_j \cdot f^{s+1-j}.
\end{displaymath}
\textcolor{black}{Ideals $I_j$ and $I_{j-1}$ are related by the following identity.}
\begin{displaymath}
	I_j = I_{j-1}\cdot f + A[s]\cdot \{\delta(a)f+(s+2-j) a \delta(f) \mid a\in I_{j-1},\delta \in\mathrm{Der}_{\mathbb C}(A)\}.
\end{displaymath}
\textcolor{black}{This formula follows from the equality} $F_jD(A) = F_1D(A) \cdot F_{j-1}D(A) = (A\oplus \mathrm{Der}_{\mathbb C}(A))\cdot F_{j-1}D(A)$. 

\vspace{0.3em}
Consequently, $0b_f(s)$ is the monic polynomial of the minimal degree such that
\begin{equation}\label{eq22}
	\textcolor{black}{b_f(s)\cdot f^{j-1} \in I_j}. 
\end{equation}
for some $j$. \textcolor{black}{This characterization will be useful in the proof of Theorem \ref{local version of Popa's conjecture}.}

\vspace{0.5em}
\textcolor{black}{Proposition \ref{unit_invariant} and Proposition \ref{local_glueing} together yield the well-defined nature of the Bernstein-Sato polynomial for non-zero effective Cartier divisors, see section \ref{introduction}.} 

\begin{proposition}\label{unit_invariant}
Suppose $X = \mathrm{Spec}\, A$ is a nonsingular \textcolor{black}{integral} complex affine variety. Let \textcolor{black}{$f\in A \setminus \{0\}$ and $u \in A$ be a unit}, then $b_f(s) = b_{uf}(s)$.
\end{proposition}
\begin{proof}
For $uf$, we define the ideals $J_j$ \textcolor{black}{analogous to} the $I_j$ above, i.e.
\begin{displaymath}
	F_jD(A)[s]  \cdot (uf)^{s+1} = J_j \cdot (uf)^{s+1-j}.
\end{displaymath}
Since $u$ is invertible, it suffices to show $I_j = J_j$ for all $j$. We \textcolor{black}{proceed} by induction. For $j = 0$, both sides \textcolor{black}{equal} $A[s]$. Now assume the identity holds for $j-1$. \textcolor{black}{The ideal $I_j$ can be expressed in terms of  $I_{j-1}$ as follows.}
\begin{displaymath}
	I_j = I_{j-1}\cdot f + A[s]\cdot \{\delta(a)f+(s+2-j) a \delta(f) \mid a\in I_{j-1},\delta \in\mathrm{Der}_{\mathbb C}(A)\}.
\end{displaymath}
The same \textcolor{black}{expression} holds for $J_j$. \textcolor{black}{That is,}
\begin{displaymath}
	J_j = J_{j-1}\cdot uf + A[s]\cdot \{\delta(a)uf+(s+2-j)a[u\delta(f)+ \delta(u)f]\mid a\in J_{j-1},\delta \in\mathrm{Der}_{\mathbb C}(A)\}.
\end{displaymath}
Since $u$ is invertible, $(s+2-j)a\delta(u) f\in J_{j-1} \cdot uf$. Then we can rewrite $J_j$ as
\begin{displaymath}
	J_j = J_{j-1}\cdot uf + A[s]\cdot \{\delta(a)uf+(s+2-j)au\delta(f)\mid a\in J_{j-1},\delta \in\mathrm{Der}_{\mathbb C}(A)\}.
\end{displaymath}
By the induction hypothesis \textcolor{black}{and the invertibility of $u$, it follows immediately that $I_j = J_j$.}
\end{proof}
\begin{remark}
\textcolor{black}{The proof of Proposition \ref{unit_invariant} is important since it illustrates a preliminary idea to use our characterization \eqref{eq22}.}
\end{remark}
\begin{remark}
\textcolor{black}{The argument in Proposition \ref{unit_invariant} also holds for Bernstein-Sato polynomials for holomorphic function germs when $f, u \in \mathbb{C}\{\bm{x}\}\setminus \{0\}$ and $u$ is invertible, thus establishing the analogous result. One can deduce from it that $ b_f(s)$ is a contact invariant of a holomorphic function germ, i.e. if $\varphi$ is an automorphism of $\mathbb C\{\bm x\}$, then $b_{f}(s) = b_{u\varphi(f)}(s)$.}
\end{remark}
\begin{proposition}[see {\cite[Proposition 2.2.6]{DMBG}}]\label{local_glueing}
Suppose $X = \mathrm{Spec}\, A$ is a nonsingular \textcolor{black}{integral complex affine} variety and $\textcolor{black}{f\in A\setminus \{0\}}$, then the global and local Bernstein-Sato polynomials are related by the formula
\begin{displaymath}
	b_f(s) = \underset{i}{\mathrm{lcm}}\, b_{f|_{U_i}}(s),
\end{displaymath}
where $\{U_i\}_i$ is an arbitrary affine open covering of $X$ and $b_{f|_{U_i}}(s)$ is the Bernstein-Sato polynomial of $f|_{U_i}$ on $U_i$.
\end{proposition}

\subsection{\textcolor{black}{Bernstein-Sato polynomials for ideals}} \label{subSec Bernstein-Sato ideal case}

The Bernstein-Sato polynomial can be generalized to an arbitrary variety\textcolor{black}{, see \cite{BMS06b,Mus22}}. The definitions in the two papers are \textcolor{black}{formulated differently but are equivalent.} 

Suppose $X = \mathrm{Spec}\, A$ is a nonsingular \textcolor{black}{integral} complex \textcolor{black}{affine} variety. Let \textcolor{black}{$f_1,...,f_r\in A$ be non-zero regular function, $\mathfrak a = (f_1,...,f_r)$ be an ideal,} and $\bm s = (s_1,...,s_r)$ be a set of \textcolor{black}{variables}. \textcolor{black}{The Bernstein-Sato polynomial of $\mathfrak a$, denoted by $b_{\mathfrak a} \in \mathbb C[\vert \bm s\vert] = \mathbb C[s_1+...+s_r]$, is the non-zero monic polynomial of the minimal degree such that}
\begin{displaymath}
b_{\mathfrak a}(\vert \bm s\vert) f_1^{s_1}f_2^{s_2} ... f_r^{s_r} \in \sum_{\bm u\in \mathbb Z^r, \vert \bm u\vert = 1} D(A)[\bm s] \cdot (\prod_{u_i < 0}\binom{s_i}{-u_i}) f_1^{s_1+u_1}f_2^{s_2+u_2}...f_r^{s_r+u_r}.
\end{displaymath}
Here, the $D(A)$-action is \textcolor{black}{natural as in \eqref{eq21}.}

The existence\textcolor{black}{, rationality,} and independence of choice of generators of $\mathfrak a$ were proved in \cite{BMS06b}. \textcolor{black}{Besides, similar to the classical Bernstein-Sato polynomial, $b_{\mathfrak a}$ serves as the algebraic avatar of the Verdier specialization complex. There is an analogue of the Kashiwara-Malgrange theorem for $b_{\mathfrak a}$, see \cite[2.8]{BMS06b}.}

\vspace{0.5em}

\textcolor{black}{The following theorem relates the newly defined object with the classical Bernstein-Sato polynomial.} 
\begin{theorem}[\cite{Mus22}]\label{bernstein-sato polynomial for ideal related to bernstein-sato polynomial}
Let $X$ be a nonsingular integral complex variety. If $f_1,...,f_r \in \Gamma(X,\mathcal O_X)$ are nonzero, generating the ideal sheaf $\mathfrak a$, and if $g = \sum_{i=1}^r f_i y_i \in \Gamma(X\times \mathbb A^r,\mathcal O_{X\times \mathbb A^r})$, then $b_{\mathfrak a} = b_g(s)/(s+1)$.
\end{theorem}
\begin{remark}
\textcolor{black}{The theorem also provides a method for computing $b_{\mathfrak a}$. However, the computation remains difficult for a general ideal $\mathfrak a$. Even when multiplicities are disregarded, determining the roots is challenging.}
\end{remark}


\vspace{0.5em}

\textcolor{black}{Next,} we review the \textcolor{black}{combinatorial descriptions} of roots of \textcolor{black}{the Bernstein-Sato polynomial of a monomial ideal. One may refer to \cite{MR2200051} and \cite{BMS06c} for more details.}

\textcolor{black}{Let $\mathfrak a \subseteq \mathbb C[\bm x]$ be a non-zero monomial ideal and} 
\begin{displaymath}
\Gamma_\mathfrak a := \{\bm u\in \mathbb N^n \mid \bm x^{\bm u} \in \mathfrak a\}
\end{displaymath}
be the semigroup associated \textcolor{black}{to} $\mathfrak a$. The Newton polyhedron of $\mathfrak a$, \textcolor{black}{denoted by} $P_{\mathfrak a}$, is defined to be the convex hull of $\Gamma_\mathfrak a$. A face $Q$ of $P_\mathfrak a$ \textcolor{black}{is defined as either} (a) $P_\mathfrak a$ itself, or (b) \textcolor{black}{a non-empty intersection of $P_{\mathfrak a}$ with a supporting affine hyperplane $H$, i.e. a hyperplane such that $P_{\mathfrak a}$
is entirely contained in one of the two closed half-spaces determined by $H$.} We call a face of type (b) a proper face.

Let $\bm e = (1,1,...,1) \in \mathbb N^n$. \textcolor{black}{For any face $Q$ of $P_{\mathfrak a}$, define}
\begin{displaymath}
M_Q:=\bm e+\mathbb N\{\bm u-\bm v \mid \bm u \in \Gamma_\mathfrak a, \bm v\in \Gamma_\mathfrak a\cap Q\}.
\end{displaymath}
\textcolor{black}{Now set} $M'_Q = \bm v_0+M_Q$, where $\bm v_0\in Q\cap \Gamma_\mathfrak a$. \textcolor{black}{Note that} $M_Q'$ does not depend on the choice of $\bm v_0$. If $Q$ is not contained in any coordinate hyperplanes, \textcolor{black}{then there exists a} linear function $L_Q$ on $\mathbb R^n$ such that $L|_Q \equiv 1$. Let $V_Q$ denote the linear span of $Q$, then $L_Q$ is uniquely determined on $V_Q$. Hence, the set 
\begin{displaymath}
R_Q :=\{-L_Q(u) \mid u\in (M_Q\setminus M_Q')\cap V_Q\}
\end{displaymath} 
is well-defined.

In particular, if $Q$ is a facet, i.e. $\text{codim}\, Q =1$, then $V_Q = \mathbb R^n$ and hence $L_Q$ is unique \textcolor{black}{on $\mathbb R^n$}. We define $m_Q$ to be the minimal positive integer such that the coefficients of $m_Q L_Q$ are all integers.

\textcolor{black}{Based on these notations, one can state the following combinatorial descriptions of roots of $b_{\mathfrak a}$ in the monomial ideal case, cf. \cite{BMS06c}.}
\begin{theorem}\label{ba monomial}
\textcolor{black}{Suppose $\mathfrak a$ is a non-zero monomial ideal. Then the set of }roots of the Bernstein-Sato polynomial of $\mathfrak a$ is \textcolor{black}{equal to} the union of \textcolor{black}{the sets $R_Q$ for the faces $Q$ of $P_{\mathfrak a}$ that are} not contained in any coordinate hyperplanes.
\end{theorem}

\begin{theorem}\label{ba monmial mod Z}
\textcolor{black}{Suppose $\mathfrak a$ is a non-zero monomial ideal. Then} the set of classes in $\mathbb Q/\mathbb Z$ of the roots of the Bernstein-Sato polynomial of $\mathfrak a$ is equal to the union of subgroups of $\mathbb Q/\mathbb Z$ generated by $\frac{1}{m_Q}$, with $Q$ running over the facets of $P_\mathfrak a$ that are not contained in any coordinate hyperplanes.
\end{theorem} 

We give \textcolor{black}{the following example to illustrate {Theorem \ref{ba monomial}}.}
\begin{example}\label{ba for monomial xaybxayb}
Let $\mathfrak a = (x^{a_1}y^{b_1},x^{a_2}y^{b_2}) \subseteq \mathbb C[x,y]$ be a monomial ideal, where $a_1 < a_2$ and $b_1 > b_2$. \textcolor{black}{The polyhedron $P_{\mathfrak a}$ has five faces: the diagonal segment $L$, the horizontal ray $X$, the vertical ray $Y$, and two points $P_1 = (a_1,b_1),P_2 = (a_2,b_2)$.} 

The linear function $L_L$ defined by $L$ is 
\begin{displaymath}
	L_L(x,y) = -\frac{b_2-b_1}{b_1a_2-a_1b_2}x+\frac{a_2-a_1}{b_1a_2-a_1b_2} y.
\end{displaymath} 
We have $M_L = (1,1)+\mathbb N^2+\mathbb Z(a_1-a_2,b_1-b_2)$. After an elementary calculation, one finds \textcolor{black}{that the elements of $M_L\setminus M_L'$ are precisely the points} \textcolor{black}{of the form} $(s+1,t+1)+\mathbb Z(a_1-a_2,b_1-b_2),s,t \geq 0$ such that
\begin{displaymath}
	\mathbb Z \cap [\frac{t-b_1}{b_2-b_1},\frac{s-a_1}{a_2-a_1}] = \varnothing.
\end{displaymath}
Consequently, $R_L$ is \textcolor{black}{given}.

For $X$ and $Y$, one finds that they contribute to $W_\mathfrak a$ the following roots.
\begin{displaymath}
	\frac{1}{a_1},...,\frac{a_1-1}{a_1},...,1,\frac{1}{b_2},...,\frac{b_2-1}{b_2},1.
\end{displaymath}

\textcolor{black}{If $a_1\neq 0$,} $L_{P_1}$ can be chosen as $\frac{x}{a_1}$. \textcolor{black}{Otherwise, $P_1$ is contained in the $y$-axis}. It is clear that $M_{P_1} \cap V_Q \cap \{x \geq a_1\} \subseteq M_{P_1}'$. So no new roots appear. \textcolor{black}{A similar argument applies to $P_2$.}

In particular, if $(a_1,b_1,a_2,b_2) = (0,b,a,0)$, we have
\begin{align*}
	L_L & = \frac{x}{a}+\frac{y}{b},\\
	M_L \setminus M_L' & = \{(s+1,t+1)\mid 0 \leq s < a, 0\leq t< b \}+\mathbb Z(-a,b),\\
	R_L & = \{-\frac{ai+bj}{ab} \mid 1\leq i\leq b,1\leq j\leq a\}.
\end{align*}
Hence, for $\mathfrak a = (x^a,y^b)$, we have 
\begin{displaymath}
	W_{\mathfrak a}  = \{-\frac{ai+bj}{ab} \mid 1\leq i\leq b,1\leq j\leq a\}.
\end{displaymath}
\end{example}

\section{The tensor property of Bernstein-Sato polynomials} \label{On Popa's Conjecture}

\subsection{\textcolor{black}{Motivations}} \label{Probable Motivations}

Let $f\in \mathbb C[\bm x]$ and $g\in \mathbb C[\bm y]$ be two non-zero polynomials. It is \textcolor{black}{clear} that $b_{f\cdot g}(s)$ \textcolor{black}{divides} $b_f(s) \cdot b_g(s)$ since if $P_f \in \mathbb C[\bm x,\bm \partial_{\bm x},s]$ and $P_g \in \mathbb C[\bm y,\bm \partial_{\bm y},s]$ \textcolor{black}{satisfy}
\begin{displaymath}
\textcolor{black}{b_f(s)\cdot f^s = P_f \cdot f^{s+1},\ b_g(s)\cdot g^s = P_g\cdot g^{s+1},}
\end{displaymath}
then we have 
\begin{displaymath}
\textcolor{black}{b_f(s)\cdot b_g(s) \cdot (f\cdot g)^s = (P_f\cdot P_g) \cdot (f\cdot g)^{s+1}.}
\end{displaymath}
So, one may be interested in the converse divisibility. In our preliminary work, we proved at least $(s+1)\mathrm{lcm}(b_f(s),b_g(s)) \mid b_{f\cdot g}(s)$. Hence, if we consider the cases when $f$ and $g$ are weighted homogeneous with isolated singularities, and all denominators of weights of $f$ \textcolor{black}{and} $g$ are coprime, we have $b_f(s)\cdot b_g(s) = b_{f\cdot g}(s)$ by {Theorem \ref{b-function for quasi-homogeneous}}.

\textcolor{black}{	The next motivation comes from the monodromy conjecture and seems more profound. Further assuming coefficients of $f$ and $g$ are rational integers, then we have
\begin{displaymath}
	Z_p(f\cdot g;s) = Z_p(f;s)\cdot Z_p(g;s).
\end{displaymath}
If the $p$-adic strong monodromy conjecture, especially the assertion about multiplicity, holds, it is natural to inquire whether $ b_f(s) \cdot b_g(s) = b_{f\cdot g}(s)$. A similar line of reasoning applies to motivic zeta functions. }

\textcolor{black}{Moreover, a similar multiplicative Thom–Sebastiani type property holds for archimedean zeta functions, see \cite{BSZ25} for definitions. Combined with \cite[Theorems 1.1 and 1.2]{BSZ25}, this yields another insight.}

\subsection{Proof of $b_{f}(s) \cdot b_g(s) = b_{f\cdot g}(s)$}\label{proof of the conjecture}


In this subsection, we give a \textcolor{black}{stronger, positive answer to  \eqref{eq2}. Specifically, we extend the setting from }polynomials $f,g$ to effective Cartier divisors on nonsingular varieties. We \textcolor{black}{start with} a simple lemma in commutative algebra.

\begin{lemma}\label{tensor product with intersection}
\textcolor{black}{Let $R$ be a principal ideal domain, and let $M$ and $N$ be free $R$-modules. Suppose $M_1, M_2 \subseteq M$ and $N_1, N_2 \subseteq N$ are submodules, which are therefore free. For any free submodules $M' \subseteq M$ and $N' \subseteq N$, there is a canonical injection $M' \otimes_R N' \to M \otimes_R N$. We may thus identify $M' \otimes_R N'$ with its image in $M \otimes_R N$. Under this identification, the following identity holds:}
\begin{displaymath}
	(M_1\cap M_2) \otimes (N_1\cap N_2) = \bigcap_{i,j\in \{1,2\}} (M_i\otimes N_j).
\end{displaymath}
\end{lemma}
\begin{proof}
\textcolor{black}{We first establish the identity:}
\begin{align}\label{eq31}
	(M_1\otimes N_1)\cap (M_2\otimes N_1) = (M_1\cap M_2) \otimes N_1.
\end{align}
Since $N_1$ is free, we \textcolor{black}{may write} $N_1 = R^{\oplus Q}$. \textcolor{black}{Fix a basis} $\{b_i\}_{i\in Q}$ \textcolor{black}{for} $N_1$, then elements of $\textcolor{black}{M_1}\otimes N_1$ can be uniquely expressed as a finite sum $\sum_{i\in Q} m_i \otimes b_i$, where $m_i\in \textcolor{black}{M_1}$. \textcolor{black}{The identity} \eqref{eq31} follows immediately.

\textcolor{black}{The same argument shows that}
\begin{displaymath}
	(M_1\otimes N_2)\cap (M_2\otimes N_2) = (M_1\cap M_2) \otimes N_2.
\end{displaymath}
Applying the same \textcolor{black}{reasoning} to $(M_1\cap M_2)$\textcolor{black}{, which is also free,} $N_1$, and $N_2$, \textcolor{black}{we conclude the proof.}
\end{proof}

The following is the local version of Theorem \ref{tensor theorem for bernstein-sato polynomial}.
\begin{theorem}\label{local version of Popa's conjecture}
Suppose $A$ and $B$ are finitely generated integral regular $\mathbb C$-algebras \textcolor{black}{and let $f\in A$ and $g\in B$ be non-zero elements.} Let $C = A\otimes_{\mathbb C} B$, and denote by $b_f(s),b_g(s)$, and $b_{f\cdot g}(s)$ the corresponding Bernstein-Sato polynomials, then \textcolor{black}{we have}
\begin{displaymath}
	b_f(s) \cdot b_g(s) = b_{f\cdot g}(s),
\end{displaymath}
\textcolor{black}{where $f\cdot g := f\otimes g \in C$.}
\end{theorem}
\begin{remark}
\textcolor{black}{We only consider the converse division, i.e. $b_f(s)\cdot b_g(s) \mid b_{f\cdot g}(s)$. Moreover, the case where either $f$ or $g$ is invertible is trivial, so we may assume $f$ and $g$ are both non-invertible.}
\end{remark}
\begin{proof}
As in {subsection \ref{subSec Bernstein-Sato polynomial, Existence and Some Results}}, we define $I_j\subseteq A[s]$ to be the ideal \textcolor{black}{satisfying}
\begin{displaymath}
	F_jD(A) \cdot f^{s+1} = I_j \cdot f^{s+1-j}.
\end{displaymath}
Similarly, \textcolor{black}{let} $J_j\subseteq B[s]$ and $L_j \subseteq C[s]$ \textcolor{black}{be the ideals defined by}
\begin{displaymath}
	F_jD(B)\cdot  g^{s+1} = J_j \cdot g^{s+1-j}\ \text{and } F_jD(C)\cdot (f\cdot g)^{s+1} = L_j\cdot (f\cdot g)^{s+1-j},
\end{displaymath}
respectively. Furthermore, we define three sequences of $\mathbb C[s]$-modules:
\begin{align*}
	M_i & := I_i \cdot f^{-i+1} \subseteq A_f[s],\\
	N_j & := J_j \cdot g^{-j+1} \subseteq B_g[s],\\
	\Lambda_l & := L_l \cdot (f\cdot g)^{-l+1} \subseteq C_{f\cdot g}[s]. 
\end{align*}
Since $I_{i} \supseteq I_{i-1}\cdot f$, \textcolor{black}{it follows that} $M_i\supseteq M_{i-1}$. \textcolor{black}{Similarly,} $N_j \supseteq N_{j-1}$ and $L_l \supseteq L_{l-1}$. By definition, we \textcolor{black}{can characterize} the Bernstein-Sato polynomials \textcolor{black}{of $f,g$, and $f\cdot g$ by}
\begin{align*}
	b_f(s) \cdot \mathbb C[s] & = \bigcup_{i = 1}^\infty (\mathbb C[s] \cap M_i),\\
	b_g(s) \cdot \mathbb C[s] & = \bigcup_{j = 1}^\infty (\mathbb C[s] \cap N_j),\\
	b_{f\cdot g}(s) \cdot \mathbb C[s] & = \bigcup_{l = 1}^\infty (\mathbb C[s] \cap \Lambda_l).
\end{align*}
It is worth mentioning that the intersections \textcolor{black}{take place in} $A_f[s]$, $B_g[s]$, and $C_{f\cdot g}[s]$\textcolor{black}{, respectively}, where $\mathbb C[s]$ can be naturally embedded. \textcolor{black}{These} identities \textcolor{black}{indicate} that for all $i,j,l \gg 0$, we have
\begin{displaymath}
	M_i\cap \mathbb C[s] = b_f(s) \cdot \mathbb C[s],\ N_j \cap \mathbb C[s] = b_g(s) \cdot \mathbb C[s],\ \Lambda_l \cap \mathbb C[s] = b_{f\cdot g}(s) \cdot \mathbb C[s].
\end{displaymath}

\textcolor{black}{Since $\mathrm{Der}_\mathbb C C = (\mathrm{Der}_\mathbb C(A) \otimes_\mathbb C B)\oplus (A \otimes_\mathbb C \mathrm{Der}_\mathbb C(B))$ and derivations of $A$ and $B$ commute, it follows that}
\begin{displaymath}
	F_lD(C) = \sum_{i+j = l} F_i D(A) \cdot F_jD(B).
\end{displaymath}
Here, we identify $F_iD(A)$ and $F_jD(B)$ with their images under canonical embeddings $F_iD_A \hookrightarrow F_i D_C$ and $F_jD_B \hookrightarrow F_j D_C$. As a result, we \textcolor{black}{obtain the following identity relating} $I_i,J_j$, and $L_l$:
\begin{displaymath}
	L_l = \sum_{i+j = l} f^{l-i} \cdot I_i \cdot g^{l-j} \cdot J_j.
\end{displaymath}
Dividing \textcolor{black}{both sides by $(f\cdot g)^{l-1}$}, we have
\begin{displaymath}
	\Lambda_l = \sum_{i+j = l} M_i \cdot N_j \subseteq M_l\cdot N_l \subseteq C_{f\cdot g}[s].
\end{displaymath}

We take $l\in \mathbb N$ sufficiently large, such that
\begin{displaymath}
	\Lambda_l \cap \mathbb C[s] = b_{f\cdot g}(s) \cdot \mathbb C[s],\ M_l\cap \mathbb C[s] = b_f(s) \cdot \mathbb C[s],\  N_l \cap \mathbb C[s] = b_g(s) \cdot \mathbb C[s].
\end{displaymath}

It suffices to prove $(M_l\cdot N_l) \cap \mathbb C[s] = (b_f(s)\cdot b_g(s)) \cdot \mathbb C[s]$. As a result, we \textcolor{black}{will} have $\Lambda_l \cap \mathbb C [s] \subseteq (b_f(s)\cdot b_g(s)) \cdot \mathbb C[s]$, i.e. $b_f(s)\cdot b_g(s)$ divides $b_{f\cdot g}(s)$. We now consider the category of $\mathbb C[s]$-modules.

The following two diagrams are Cartesian products\textcolor{black}{, with all eight arrows being monomorphisms of $\mathbb C[s]$-modules.}
\begin{displaymath}
	\xymatrix{
		b_f(s)\cdot \mathbb C[s] \ar[r] \ar[d] & \mathbb C[s] \ar[d]\\
		M_l \ar[r] & A_f[s],
	}
\end{displaymath}
\begin{displaymath}	
	\xymatrix{
		b_g(s)\cdot \mathbb C[s] \ar[r] \ar[d] & \mathbb C[s] \ar[d]\\
		N_l \ar[r] & B_g[s].
	}
\end{displaymath}
Tensoring the two squares \textcolor{black}{over $\mathbb C[s]$ yields} the following commutative diagram.
\begin{align}\label{eq311}
	\xymatrix{
		(b_f(s)\cdot \mathbb C[s]) \otimes_{\mathbb C[s]} (b_g(s)\cdot \mathbb C[s]) \ar[rr] \ar[d] && \mathbb C[s] \otimes_{\mathbb C[s] } \mathbb C[s] \ar[d]\\
		M_l \otimes_{\mathbb C[s]} N_l \ar[rr] && A_f[s] \otimes_{\mathbb C[s]} B_g[s].
	}
\end{align}
A simple observation shows $M_l,N_l,A_f[s]$, and $B_g[s]$ are all $\mathbb C[s]$-free modules. \textcolor{black}{Indeed}, $A_f[s] = A_f \otimes_{\mathbb C} \mathbb C[s]$ is \textcolor{black}{a tensor product of a vector space with $\mathbb C[s]$}, and \textcolor{black}{is} hence free. Moreover, since $\mathbb C[s]$ is a principal ideal domain and $M_l$ is a $\mathbb C[s]$-submodule of $A_f[s]$, $M_l$ is also free. \textcolor{black}{The same argument applies to $N_l$ and $B_g[s]$.} Consequently, one deduces that the four arrows in \eqref{eq311} are all \textcolor{black}{injective}.

Under the canonical isomorphism $A_f[s] \otimes_{\mathbb C[s]} B_g[s] \simeq C_{f\cdot g} [s]$, the diagram \eqref{eq311} is identified with the following:
\begin{displaymath}
	\xymatrix{
		(b_f(s)\cdot b_g(s))\cdot \mathbb C[s] \ar[r] \ar[d] & \mathbb C[s] \ar[d]\\
		M_l \cdot N_l \ar[r] & C_{f\cdot g}[s].
	}
\end{displaymath}
Therefore, we only need to show \eqref{eq311} is a Cartesian product.

Applying {Lemma \ref{tensor product with intersection}}, we have
\begin{displaymath}
	(b_f(s)\cdot \mathbb C[s]) \otimes_{\mathbb C[s]} (b_g(s)\cdot \mathbb C[s]) = (M_l \otimes_{\mathbb C[s]} N_l) \cap (\mathbb C[s] \otimes_{\mathbb C[s] } \mathbb C[s]) \cap (M_l \otimes_{\mathbb C[s]} \mathbb C[s]) \cap (\mathbb C[s] \otimes _{\mathbb C[s]} N_l).
\end{displaymath}
Since all arrows in \eqref{eq311} are injective, it suffices to show
\begin{align}\label{eq32}
	(M_l \otimes_{\mathbb C[s]} N_l) \cap (\mathbb C[s] \otimes_{\mathbb C[s] } \mathbb C[s]) \subseteq M_l \otimes_{\mathbb C[s]} \mathbb C[s]. 
\end{align}
Similarly, the corresponding \textcolor{black}{containment} holds for $N_l$. As a result, \eqref{eq311} is a Cartesian product. 

To prove \eqref{eq32}, we use again the canonical isomorphism $A_f[s] \otimes_{\mathbb C[s]} B_g[s] \simeq C_{f\cdot g} [s]$ to pass everything to $C_{f\cdot g}[s]$. It \textcolor{black}{suffices to show that}
\begin{equation}\label{eq33k}
	(M_l \cdot N_l) \cap \mathbb C[s] \subseteq M_l.
\end{equation}

Let $\sum_i a_i \cdot f^{-l+1} \cdot b_i \cdot g^{-l+1} = c(s)$ be an element in the left-hand side of \eqref{eq33k}, where $a_i \in I_l, b_i \in J_l$, and $c(s)\in \mathbb C[s]$. That is,
\begin{align}\label{eq33}
	\sum_i a_ib_i = c(s) (f\cdot g)^{l-1}. 
\end{align}
The equality \eqref{eq33} \textcolor{black}{takes place in} $ C[s]$.

We \textcolor{black}{now} localize $B$ to show \textcolor{black}{that} $c(s)\in M_l$. Let $\mathfrak m_{\bm y}\subseteq B$ be a maximal ideal such that $g\in \mathfrak m_{\bm y}$, \textcolor{black}{where} $\bm y = (y_1,...,y_m)$ \textcolor{black}{denotes} a \textcolor{black}{system} of local coordinates of $B$ at $\mathfrak m_{\bm y}$. Let $\tilde B = B_{\mathfrak m_{\bm y}}$ be the localization, and \textcolor{black}{denote its maximal ideal again by} $\mathfrak m_{\bm y}$. Tensoring \textcolor{black}{the natural embedding $B \hookrightarrow \tilde B$ with $A$ and $\mathbb C[s]$}, we have $C[s] \hookrightarrow  (A \otimes_{\mathbb C} \tilde B)[s]$.
Suppose the order of $g$ in $\tilde B$ is $t > 0$, i.e. $g\in \mathfrak m_{\bm y}^t$ but $g\not \in \mathfrak m_{\bm y}^{t+1}$. Since $B$ is a regular local ring, \textcolor{black}{the associated graded ring}
\begin{align*}
	\mathrm{gr}_{\mathfrak m_{\bm y}} \tilde B = \bigoplus_{j\geq 0} \mathfrak m_{\bm y}^j/\mathfrak m_{\bm y}^{j+1}. 
\end{align*}
is isomorphic to \textcolor{black}{a} polynomial ring \textcolor{black}{in} $m$ variables, with variables \textcolor{black}{corresponding to the} images of $\bm y$ in $\mathrm{gr}_{\mathfrak m_{\bm y}} \tilde B$. Suppose $\bm y^{\bm \alpha}$ appears in the image of \textcolor{black}{$g$ in $\tilde B/\mathfrak m_{\bm y}^{t+1}$} and \textcolor{black}{assume, without loss of generality, that its} coefficient is $1$.

Let us \textcolor{black}{consider the equation} \eqref{eq33} in $(A\otimes_{\mathbb C} \tilde B)[s]$. \textcolor{black}{Reducing modulo} $\mathfrak m_{\bm y}^{(l-1)t+1}$, \textcolor{black}{the} identity also holds in $A[s] \otimes_{\mathbb C[s]} (\tilde B/\mathfrak m_{\bm y}^{(l-1)t+1}) [s]$. Since $(\tilde B/\mathfrak m_{\bm y}^{(l-1)t+1}) [s]$ is a \textcolor{black}{free} $\mathbb C[s]$-module, we \textcolor{black}{may} compare the coefficients of $\bm y^{(l-1)\bm \alpha}$, which are in $A[s]$, on both sides of \textcolor{black}{\eqref{eq33}}. On \textcolor{black}{the right-hand side}, the coefficient of $\bm y^{(l-1)\bm \alpha}$ is $c(s) f^{l-1}$. On \textcolor{black}{the left-hand side}, \textcolor{black}{let $\tilde b_i(s)\in \mathbb C[s]$ be the coefficient of $\bm y^{(l-1)\bm \alpha}$ in $b_i$, then} we obtain $\sum_i \tilde b_i(s) a_i = c(s) f^{l-1}$. \textcolor{black}{Therefore}, $c(s) = f^{-l+1}\cdot (\sum_i\tilde b_i(s) a_i) \in M_l$\textcolor{black}{, which completes the proof.}
\end{proof}

As \textcolor{black}{a} corollary, we can generalize \textcolor{black}{Theorem \ref{local version of Popa's conjecture}} to a global one. \textcolor{black}{To do this,} we first define the tensor \textcolor{black}{product} of Cartier divisors, as a generalization of $f\cdot g$.
\begin{definition}\label{tensor of divisors}
Suppose $X_1,X_2$ are nonsingular complex varieties and $D_1,D_2$ are Cartier divisors on $X_1,X_2$ respectively. We define \textcolor{black}{the Cartier divisor $D_1\otimes D_2$ on $X_1\times X_2$ in the following way. Choose} affine coverings $\{U_i\},\{V_j\}$ of $X_1$ and $X_2$ respectively, such that \textcolor{black}{$D_1|_{U_i} = \mathrm{div}(f_i)$ and $D_2|_{V_j} = \mathrm{div}(g_j)$. Then $D_1\otimes D_2$ is defined locally on each $U_i\times V_j$ by $f_i\otimes g_j$.}
\end{definition}
\begin{proof}[Proof of Theorem \ref{tensor theorem for bernstein-sato polynomial}]
As in {Definition \ref{tensor of divisors}}, \textcolor{black}{we} take open affine coverings $\{U_i\},\{V_j\}$. By {Theorem \ref{local version of Popa's conjecture}}, we have
\begin{displaymath}
	b_{D_1\otimes D_2|_{U_i\times V_j}}(s) = b_{D_1|_{U_i}}(s) \cdot b_{D_2|_{V_j}}(s).
\end{displaymath}
Since $X_1 \times V_j = \bigcup_i U_i\times V_j$, we have
\begin{displaymath}
	b_{D_1\otimes D_2|_{X_1\times V_j}}(s) = \underset{i}{\mathrm{lcm}}\, b_{D_1|_{U_i}}(s)\cdot b_{D_2|_{V_j}}(s) = [\underset{i}{\mathrm{lcm}}\, b_{D_1|_{U_i}}(s)]\cdot b_{D_2|_{V_j}}(s) = b_{D_1}(s) \cdot  b_{D_2|_{V_j}}(s).
\end{displaymath}
Furthermore, by \textcolor{black}{the covering} $X_1 \times X_2 = \bigcup_j X_1 \times V_j$, we have
\begin{displaymath}
	b_{D_1\otimes D_2}(s) = \underset{j}{\mathrm{lcm}}\, b_{D_1}(s)\cdot b_{D_2|_{V_j}}(s) = b_{D_1}(s)\cdot [\underset{j}{\mathrm{lcm}}\, b_{D_2|_{V_j}}(s)] = b_{D_1}(s) \cdot  b_{D_2}(s).
\end{displaymath}
\end{proof}

\subsection{Generalization to ideals}\label{Generalization to Ideals}

Having solved the problem about the product of two regular functions, \textcolor{black}{a natural generalization arises as follows.}
\begin{question}\label{TS:for:ideals}
For \textcolor{black}{non-zero ideals} $\mathfrak a \subseteq \mathbb C[\bm x]$ and $\mathfrak b\subseteq \mathbb C[\bm y]$, is $b_{\mathfrak a\cdot \mathfrak b} = b_{\mathfrak a}\cdot b_{\mathfrak b}$?
\end{question} 


\textcolor{black}{If either $\mathfrak a$ or $\mathfrak b$ is principal, the expected identity in Question \ref{TS:for:ideals} holds, see Proposition \ref{generalized Popa's conjecture in a case}. Here is a proof.}
\begin{proof}[Proof of Proposition \ref{generalized Popa's conjecture in a case}]
\textcolor{black}{Suppose $f_1,...,f_r$ are non-zero generators of $\mathfrak a$.} By {Theorem \ref{bernstein-sato polynomial for ideal related to bernstein-sato polynomial}}, $b_\mathfrak a = b_{h}(s)/(s+1)$, where $h = \sum f_i z_i$ and \textcolor{black}{$z_1,...,z_r$} are \textcolor{black}{new variables}. Again by {Theorem \ref{bernstein-sato polynomial for ideal related to bernstein-sato polynomial}}, $b_{\mathfrak a \cdot (g)} = b_{h'}(s)/(s+1)$, where $h' = \sum f_i g z_i$. Applying {Theorem \ref{local version of Popa's conjecture}}, we have $b_{h'} = b_g \cdot b_h$, and hence $b_{\mathfrak a\cdot (g)} = b_\mathfrak a \cdot b_g$.
\end{proof}

\textcolor{black}{The question is still interesting even \textcolor{black}{if} we neglect the multiplicities of roots.} Recall that $W_\mathfrak n$ denotes the roots of $b_\mathfrak n$ for an ideal $\mathfrak n$ contained in \textcolor{black}{the} polynomial ring. \textcolor{black}{In the remainder of this subsection, we will prove Theorem \ref{bab for monomial ideals1}.}


\vspace{0.5em}

\textcolor{black}{To facilitate our subsequent discussion, we first establish several properties of the Newton polyhedron associated with the product of two monomial ideals.} 


\begin{lemma}\label{Newton Diagram of product monomial ideal}
\textcolor{black}{Let $\mathfrak a \subseteq \mathbb C[\bm x]$ and $\mathfrak b\subseteq \mathbb C[\bm y]$ be non-zero monomial ideals. Denote their associated semigroups by $\Gamma_{\mathfrak a}$ and $\Gamma_{\mathfrak b}$, and their Newton polyhedra by $P_{\mathfrak a}$ and $P_{\mathfrak b}$. Then: }

\noindent(a) $\Gamma_{\mathfrak a\cdot \mathfrak b} = \Gamma_\mathfrak a\times \Gamma_{\mathfrak b}$.

\noindent (b) $P_{\mathfrak{a}\cdot \mathfrak{b}} = P_\mathfrak a \times P_\mathfrak b $.

\noindent (c) Faces of $P_{\mathfrak a \cdot \mathfrak b}$ are exactly those of \textcolor{black}{the form} $Q_1 \times Q_2$, where $Q_1$(resp. $Q_2$) is a face of $P_\mathfrak a$ (resp. $P_\mathfrak b$).

\noindent (d) For \textcolor{black}{faces} $Q_1 \subseteq P_\mathfrak a$ and $Q_2 \subseteq P_\mathfrak b$, we have
\begin{displaymath}
M_{Q_1\times Q_2}\setminus M'_{Q_1\times Q_2}  = [(M_{Q_1}\setminus M_{Q_1}')\times M_{Q_2}] \cup [M_{Q_1} \times (M_{Q_2}\setminus M_{Q_2}')].
\end{displaymath}
\end{lemma}
\begin{proof}
We prove only part (c), \textcolor{black}{as (a) and (b) are straightforward and (d) follows from a direct computation on sets.}

It suffices to consider proper faces. \textcolor{black}{Let $H$ be a hyperplane giving a proper face of $P_{\mathfrak a \cdot \mathfrak b}$.} Suppose $H$ is given by the equation: 
\begin{displaymath}
\lambda_1x_1+...\lambda_nx_n+\mu_1 y_1+...+\mu_m y_m - C = 0.
\end{displaymath}

\textcolor{black}{One can show} that all $\lambda_i,\mu_j,C\in \mathbb Q$ have the same sign \textcolor{black}{by using the existence of a positive integer $\theta$  such that $(\theta,\theta,...,\theta)+\mathbb N^{n+m} \subseteq P_{\mathfrak{a}\cdot \mathfrak{b}}$.} 
\textcolor{black}{Let $C_1$ and $C_2$ be the minimum of $X = \lambda_1x_1+...\lambda_nx_n$ and $Y = \mu_1 y_1+...+\mu_m y_m$ on $P_{\mathfrak a}$ and $P_{\mathfrak b}$ respectively. Then} we have $C = C_1+C_2$ and $H\cap P_{\mathfrak a \cdot \mathfrak b} = H_1 \times H_2$, where $H_1$ (resp. $H_2$) is the face \textcolor{black}{of $P_{\mathfrak a}$ (resp. $P_{\mathfrak b}$) given} by $X-C_1$ (resp. $Y-C_2$).
\end{proof}

\textcolor{black}{Now we can prove {Theorem \ref{bab for monomial ideals1}}.} 
\begin{proof}[Proof of Theorem \ref{bab for monomial ideals1}]

\textcolor{black}{By Lemma \ref{Newton Diagram of product monomial ideal}(c), there are two types of proper faces of $P_{\mathfrak a\cdot \mathfrak b}$ not contained in any coordinate hyperplanes: (i) $Q_1 \times P_{\mathfrak b}$ or $P_{\mathfrak a} \times Q_2$, and (ii) $Q_1 \times Q_2$, where $Q_1$ (resp. $Q_2$) is a proper face of $P_{\mathfrak a}$ (resp. $P_{\mathfrak b}$) not contained in any coordinate hyperplanes. }

\textcolor{black}{\textit{Type (i).} Let $Q = Q_1$ and $L_Q$ be a linear function on $\mathbb R^n$ such that $L_Q|_Q \equiv 1$, then $L_Q\oplus 0$ is a linear function which takes value $1$ on $Q\times P_\mathfrak b$. By {Lemma \ref{Newton Diagram of product monomial ideal}(d)} and the fact that $M_{P_\mathfrak b} = \mathbb Z^m$, we have $M_{Q\times P_\mathfrak b} = (M_Q\setminus M_Q')\times \mathbb Z^m$. Since the linear span of $Q\times P_\mathfrak b$ is $V_Q\times \mathbb R^m$, it follows that $R_{Q\times P_\mathfrak b} = R_Q$. The same argument works for $Q_2$. Hence, by Theorem \ref{ba monomial}, the \textcolor{black}{subset} of roots of $b_{\mathfrak a\cdot \mathfrak b}$ coming from faces of type (i) is exactly $W_{\mathfrak a}\cup W_{\mathfrak b}$.}

\textcolor{black}{\textit{Type (ii).}  By Lemma \ref{Newton Diagram of product monomial ideal}(d), we have $M_{Q_1\times Q_2}\setminus M'_{Q_1\times Q_2}  = [(M_{Q_1}\setminus M_{Q_1}')\times M_{Q_2}] \cup [M_{Q_1} \times (M_{Q_2}\setminus M_{Q_2}')]$. We claim that all roots originating from $(M_{Q_1}\setminus M_{Q_1}')\times M_{Q_2}$ lie in $W_{\mathfrak a}$. To see this, choose $L_{Q_1\times Q_2} = L_{Q_1} \oplus 0$, note that $V_{Q_1 \times Q_2} \subseteq V_{Q_1} \times V_{Q_2}$, and use Theorem \ref{ba monomial}. By symmetry, no new roots arise from $M_{Q_1} \times (M_{Q_2}\setminus M_{Q_2}')$, either.}

\textcolor{black}{We conclude the proof by combining the results for two types. }
\end{proof}

\begin{remark}
\textcolor{black}{If one considers merely classes in $\mathbb Q/\mathbb Z$ of the roots, the expected identity will be easier to prove, simply by applying Theorem \ref{ba monmial mod Z} and noticing that all facets of $\mathbb P_{\mathfrak a\cdot \mathfrak b}$ not contained in any coordinate hyperplanes are exactly those of the type (i) in the proof of Theorem \ref{bab for monomial ideals1}.}
\end{remark}

\begin{thebibliography}{AMHJ+22}

\bibitem[AC75]{AC} N. A’Campo. La fonction z\^eta d’une monodromie. {\em Comment. Math. Helv.,} 50: 233-248, 1975.

\bibitem[AMHJ$^+$22]{doi:10.1142/S0219199721500838}
J.~\`{A}lvarez Montaner, D.~Hern\'{a}ndez, J.~Jeffries, L.~N\'{u}\~{n}ez
Betancourt, P.~Teixeira, and E.~E. Witt.
\newblock Bernstein–Sato functional equations, V-filtrations, and multiplier
ideals of direct summands.
\newblock {\em Commun. Contemp. Math.}, 24(10): 2150083,
2022.


\bibitem[Ber71]{MR0290097}
I.~N. Bern\v{s}te\u{\i}n.
\newblock Modules over a ring of differential operators. {A}n investigation of
the fundamental solutions of equations with constant coefficients.
\newblock {\em Funkcional. Anal. i Prilo\v{z}en.}, 5(2): 1--16, 1971.

\bibitem[Bla24]{Bla24}
G. Blanco.
\newblock Topological roots of the Bernstein-Sato polynomial of plane curves.
\newblock{\em \url{https://arxiv.org/abs/2406.09034}}, 2024.

\bibitem[Bud12]{Bud12}
N. Budur.
\newblock Singularity invariants related to {M}ilnor fibers: survey.
\newblock In {\em Zeta functions in algebra and geometry}, volume 566 of {\em
	Contemp. Math.}, pages 161--187. Amer. Math. Soc., Providence, RI, 2012.

\bibitem[BMS06a]{MR2200051}
N.~Budur, M.~Musta\c{t}\u{a}, and M.~Saito.
\newblock Roots of {B}ernstein-{S}ato polynomials for monomial ideals: a
positive characteristic approach.
\newblock {\em Math. Res. Lett.}, 13(1): 125--142, 2006.

\bibitem[BMS06b]{BMS06b}
N.~Budur, M.~Musta\c{t}\v{a}, and M.~Saito.
\newblock Bernstein-{S}ato polynomials of arbitrary varieties.
\newblock {\em Compos. Math.}, 142(3): 779--797, 2006.

\bibitem[BMS06c]{BMS06c}
N.~Budur, M.~Musta\c{t}\v{a}, and M.~Saito.
\newblock Combinatorial description of the roots of the {B}ernstein-{S}ato
polynomials for monomial ideals.
\newblock {\em Comm. Algebra}, 34(11): 4103--4117, 2006.


\bibitem[BSZ25]{BSZ25}
N.~Budur, Q.~Shi, and H.~Zuo.
\newblock Polar loci of multivariable archimedean zeta functions.
\newblock{\em \url{https://arxiv.org/abs/2504.10051}}, 2025.

\bibitem[BVWZ21]{Zero_Loci_of_Bernstein-Sato_Ideals}
N. Budur, R. van~der Veer, L. Wu, and P. Zhou.
\newblock Zero loci of {B}ernstein-{S}ato ideals.
\newblock {\em Invent. Math.}, 225: 45-72, 2021.

\bibitem[CDNM24]{CDNM}
A. Castaño Domínguez and L. Narváez Macarro.
\newblock{ On the reduced Bernstein-Sato polynomial of Thom-Sebastiani singularities}.
\newblock {\em Rev. Mat. Complut.}, 37: 877–886, 2024.

\bibitem[DLY24]{DLY24}
D. Davis, A. Lőrincz, and R. Yang. 
\newblock{Archimedean zeta functions, singularities, and Hodge theory}.
\newblock{\em \url{https://arxiv.org/abs/2412.07849v1}}.

\bibitem[DL98]{DL98}
J. Denef and F. Loeser.
\newblock Motivic {I}gusa zeta functions.
\newblock {\em J. Algebraic Geom.}, 7: 505-537, 1998.


\bibitem[ELSV04]{ELSV04}
L.~Ein, R.~Lazarsfeld, K.~E. Smith, and D.~Varolin.
\newblock Jumping coefficients of multiplier ideals.
\newblock {\em Duke Math. J.}, 123(3): 469--506, 2004.


\bibitem[Gra10]{MR2766095}
M.~Granger.
\newblock Bernstein-{S}ato polynomials and functional equations.
\newblock In {\em Algebraic approach to differential equations}, pages
225--291. World Sci. Publ., Hackensack, NJ, 2010.


\bibitem[Igu75]{Igu75}
J.~Igusa.
\newblock Complex powers and asymptotic expansions. {II}. {A}symptotic
expansions.
\newblock {\em J. Reine Angew. Math.}, 278/279: 307--321, 1975.

\bibitem[Igu00]{Igu00}
J.~Igusa.
\newblock An introduction to the theory of local zeta functions.
\newblock{\em AMS/IP Studies in Advanced Mathe
	matics., 14. American Mathematical Society, Providence, RI; International Press, Cambridge, MA}, 2000.


\bibitem[Kas83]{Kashiwara_Malgrange_Theorem_Kashiwara}
M.~Kashiwara.
\newblock Vanishing cycle sheaves and holonomic systems of differential
equations.
\newblock {\em Algebraic geometry, {Proc}. {Jap}.-{Fr}. {Conf}., {Tokyo} and {Kyoto}
	1982}, {Lect}. {Notes} {Math}., 1016: 134-142, 1983.


\bibitem[KS11]{KS11}
S. Kov\'{a}cs and K. Schwede.
\newblock Hodge theory meets the minimal model program: a survey of log
canonical and {D}u {B}ois singularities.
\newblock In {\em Topology of stratified spaces}, volume~58 of {\em Math. Sci.
	Res. Inst. Publ.}, pages 51--94. Cambridge Univ. Press, Cambridge, 2011.


\bibitem[Lee24]{lee2024multiplicative}
J. Lee.
\newblock Multiplicative Thom-Sebastiani for Bernstein-Sato polynomials.
\newblock{\em \url{https://arxiv.org/abs/2402.04512}}, 2024.


\bibitem[Mal83]{Kashiwara_Malgrange_Theorem_Malgrange}
B. Malgrange.
\newblock Polyn{\^o}mes de {Bernstein}-{Sato} et cohomologie evanescente.
\newblock {\em Ast{\'e}risque} 101-102: 243-267, 1983.


\bibitem[Mus22]{Mus22}
M.~Musta\c{t}\u{a}.
\newblock Bernstein-{S}ato polynomials for general ideals vs. principal ideals.
\newblock {\em Proc. Amer. Math. Soc.}, 150(9): 3655--3662, 2022.

\bibitem[Nic10]{Nic10}
J.~Nicaise.
\newblock An introduction to {$p$}-adic and motivic zeta functions and the
monodromy conjecture.
\newblock In {\em Algebraic and analytic aspects of zeta functions and
	{$L$}-functions}, volume~21 of {\em MSJ Mem.}, pages 141--166. Math. Soc.
Japan, Tokyo, 2010.

\bibitem[Pop18]{DMBG}
M.~Popa.
\newblock {$\mathscr D$}-modules in birational geometry.
\newblock {\em \url{https://people.math.harvard.edu/~mpopa/notes/DMBG-posted.pdf},
	2018.}

\bibitem[Sai93]{Sai93}
M.~Saito.
\newblock On {$b$}-function, spectrum and rational singularity.
\newblock {\em Math. Ann.}, 295(1): 51--74, 1993.

\bibitem[Sai06]{saito2006introduction}
M.~Saito.
\newblock Introduction to a theory of {$b$}-functions.
\newblock{\em \url{https://arxiv.org/abs/math/0610783v1}}, 2006.

\bibitem[Sai09]{Sai09}
M.~Saito.
\newblock On the {H}odge filtration of {H}odge modules.
\newblock {\em Mosc. Math. J.}, 9(1): 161--191, back matter, 2009.

\bibitem[Sai16]{MR3573952}
M.~Saito.
\newblock Bernstein-{S}ato polynomials of hyperplane arrangements.
\newblock {\em Selecta Math. (N.S.)}, 22(4): 2017--2057, 2016.



\bibitem[SY23]{SY23}
C. Schnell and R. Yang.
\newblock Higher multiplier ideals.
\newblock{\em \url{https://arxiv.org/abs/2309.16763}}, 2023.


\bibitem[Vey25]{Vey25}
			W. Veys.
			\newblock Introduction to the monodromy conjecture.
            \newblock{\em Handbook of geometry and topology of singularities VII}, Springer, Cham, 721–765, 2025.

\bibitem[Wal15]{Walter_PhD}
R.~Walters.
\newblock {An application of {D}-modules: {T}he {B}ernstein-{S}ato
	polynomial}.
\newblock{\em ProQuest LLC, Ann Arbor, MI}, 2015.

\bibitem[Wei67]{MR0234930}
A.~Weil.
\newblock {Basic number theory}. 
\newblock{\em volume Band 144 of {\em Die Grundlehren
		der mathematischen Wissenschaften},} Springer-Verlag New York, Inc., New York, 1967.



\end{thebibliography}
\end{document}